\documentclass[11pt]{amsart}
\usepackage{amsmath}
\newtheorem{thm}{Theorem}[section]

\newtheorem{lem}[thm]{Lemma}
\newtheorem{prop}[thm]{Proposition}
\theoremstyle{definition}
\newtheorem{defn}[thm]{Definition}
\theoremstyle{remark}
\newtheorem{rem}[thm]{\bf Remark}

\numberwithin{equation}{section}
\begin{document}
\title[Reversibility in the groups $\textrm{PL}^+(\mathbb{S}^{1})$ and $\textrm{PL}(\mathbb{S}^{1})$
%GROUP OF PIECEWISE LINEAR HOMEOMORPHISMS
]
%{REVERSIBILITY IN THE GROUP OF PIECEWISE LINEAR HOMEOMORPHISMS OF THE CIRCLE}
{Reversibility in the groups $\textrm{PL}^+(\mathbb{S}^{1})$ and
$\textrm{PL}(\mathbb{S}^{1})$}

\author{ Khadija Ben Rejeb and Habib Marzougui %$^{\dag}$
}
%\footnote{$^{\dag}$ Senior Associate of ICTP.}
\address{ Habib Marzougui, \
University of Carthage, Faculty of Science of Bizerte, Department of Mathematics, Jarzouna. 7021. Tunisia.
E.mail:hmarzoug@ictp.it; habib.marzougui@fsb.rnu.tn}
\address{ Khadija Ben Rejeb,\  Institut Sup\'erieur d'Informatique et des Technologies de Communication de Hammam Sousse, \
Route principale 1, 4011 Hammam Sousse, Tunisia. \\
E.mail: kbrjeb@yahoo.fr}



\begin{abstract}
\noindent  Let $\textrm{PL}^+(\mathbb{S}^{1})$ be the group of order
preserving piecewise linear homeomorphisms of the circle. An element
in $\textrm{PL}^+(\mathbb{S}^{1})$ is called {\it reversible} in
$\textrm{PL}^+(\mathbb{S}^{1})$ if it is conjugate to its inverse in
$\textrm{PL}^+(\mathbb{S}^{1})$. We characterize the reversible
elements in $\textrm{PL}^+(\mathbb{S}^{1})$. We also perform a
similar characterisation in the full group
$\textrm{PL}(\mathbb{S}^{1})$ of piecewise linear homeomorphisms of
the circle.
\end{abstract}
\maketitle

\section{\bf Introduction}
 Let $G$ be a group.  An element $g \in G$ is called {\it reversible} in $G$ if it is conjugate to its inverse in $G$;
there exists $h \in G$ such that  \begin{equation} g^{-1} = hgh^{-1}. \end{equation}
 We say that $g$ is reversed by $h$. An element $h$ in $G$ is called an \emph{involution} if $h = h^{-1}$.
 If the conjugating element $h$ can be chosen to be an involution then $g$
is called {\it strongly reversible}. Any strongly reversible element
can be expressed as a product of two involutions. So involutions are
strongly reversible and strongly reversible elements are reversible.
For $n\in \mathbb{N}^{*}$, define $$I_n(G) =
 \{\tau_1 \tau_2\dots  \tau_n \ \ :  \ \ \forall \ 1 \leq i \leq n, \ \tau_i \ \hbox{is an involution in} \ G \};$$
$$R_n(G) = \{\ g_1 g_2\dots g_{n}  \ \ :  \ \ \forall \ 1 \leq i \leq n, \ g_i \
\hbox{is a reversible element in} \ G\}.$$ \noindent For each
integer $n \in \mathbb{N}^{*}$, it is clear that $I_n(G) \subseteq
I_{n+1}(G)$ and $R_n(G) \subseteq R_{n+1}(G)$. The set $I_1(G)$
(resp. $I_2(G)$) consists of the involutions (resp. the strongly
reversible elements) in $G$. Denote by

\textbullet \ $\mathbb{S}^{1} = \{z \in \mathbb{C} \ | \ \left|z\right| = 1\}$ the circle, it is a multiplicative group.

\textbullet \  $\textrm{Homeo}(\mathbb{S}^{1})$ (resp.
$\textrm{Homeo}^+(\mathbb{S}^1)$) the group of all homeomorphisms
(resp. orientation-preserving homeomorphisms) of $\mathbb{S}^1$.

 \textbullet \ $\textrm{id}$ the identity map of $\mathbb{S}^1$.
 
\textbullet \ $s:  \mathbb{S}^{1}\longrightarrow \mathbb{S}^{1}, \  z\mapsto \overline{z}$ the reflection.

\textbullet \ $\textrm{Fix}(f)$ the set of fixed points of $f$.

 \noindent Two elements $f$ and $g$ of $\textrm{Homeo}(\mathbb{S}^1)$ are called \textit{conjugate }in $\textrm{Homeo}(\mathbb{S}^1)$ if 
there exists $h \in
\textrm{Homeo}(\mathbb{S}^1)$ such that $g = hfh^{-1}$. It is well
known that for every homeomorphism $f: \mathbb{S}^{1}\rightarrow
\mathbb{S}^1$ there exists a unique
 (up to a translation by an integer) homeomorphism $\tilde{f} : \mathbb{R} \rightarrow \mathbb{R}$ such that
\begin{equation}f(e^{i2 \pi x}) = e^{i2 \pi \tilde{f}(x)} \noindent \\ \textrm{ and} \ \tilde{f}(x+1) = \tilde{f}(x)+k,
\textrm{ for all } x \in \mathbb{R}\end{equation} where $k \in \{-1,
1\}$. Such a homeomorphism $\tilde{f}$ is called a \emph{lift} of
$f$. We call $f$ orientation-preserving if $k = 1$; resp.
orientation-reversing if $k = -1$; which is equivalent to the fact
that $\tilde{f}$ is increasing, resp. $\tilde{f}$ is decreasing.
\bigskip
\medskip

 \begin{defn}\label{d1} A homeomorphism $f$ of $\mathbb{S}^{1}$ is said to be {\it piecewise
linear} ($\mathrm{PL}$) if it is 
derivable except at finitely or countably many points  $(c_{i})_{i\in \mathbb{N}}$ called
\textit{break points} of $f$  at which $f$ admits left and right
derivatives (denoted, respectively, by  Df$_{-}$  and  Df$_{+}$)
and such that the derivative  $\mathrm{Df}: \ S^{1} \longrightarrow
\mathbb{R}^{*}$ is constant on each connected component of $\mathbb{S}^{1}\backslash  \{c_{i}: \ i\in \mathbb{N}\}$.
\end{defn}
\medskip

Let $f$ be a {\it piecewise linear} ($\textrm{PL}$) homeomorphism
$f$ of $\mathbb{S}^{1}$. Define $\sigma_f(x) = {\textrm{Df}_-(x)
\over \textrm{Df}_+(x)}$ the $f$-jump in $x\in \mathbb{S}^1$. Denote by $B(f) = \{c_{i}: \ i\in \mathbb{N}\}$ the set of 
\emph{break points} of $f$. So $B(f)= \{x\in \mathbb{S}^1: \sigma_f(x)\neq
 1\}$.\\

%\begin{defn}\label{d2} **** A homeomorphism $f$ of $\mathbb{S}^{1}$ is said to be {\it piecewise linear with finite break points}
%($\textrm{PLF}$) provided that $f$ is $\textrm{PL}$ with $B(f)$
%finite.
%\end{defn}
%\medskip

 %\begin{defn} \label{d22} A homeomorphism $f$ of $\mathbb{S}^{1}$ is said to be {\it piecewise linear} ($\textrm{PL}$) if there exists
 %a lift $\tilde{f}$ of $f$ to $\mathbb{R}$ and a finite subdivision
%0<a_1<a_2< \dots < a_p=1$ of the interval $[0,1]$ such that
%\begin{equation}\tilde{f}_{|[a_i, a_{i+1}]}(x) = \lambda_i x +
%\beta_i, \ \ \lambda_i, \beta_i \in \mathbb{R}. \end{equation}
%\end{defn}
\medskip

 \noindent The homeomorphism $f$ is $\textrm{PL}$ (resp. $\textrm{PL}^+$, $\textrm{PL}^-$) if and only if  $\tilde{f}$
 is a piecewise linear (resp. piecewise linear increasing, piecewise linear decreasing) homeomorphism of the real
line $\mathbb{R}$.

Denote by

\textbullet \  \  $\textrm{PL}(\mathbb{S}^{1})$ the group of all
piecewise linear homeomorphisms of $\mathbb{S}^1$,

\textbullet \ $\textrm{PL}^+(\mathbb{S}^{1})$ the group of orientation-preserving elements of $\textrm{PL}(\mathbb{S}^{1})$,

\textbullet \  \ $\textrm{PL}^-(\mathbb{S}^{1})$ the set of
orientation-reversing elements of $\textrm{PL}(\mathbb{S}^{1})$.

%\textbf{****} \textbullet \  \ $\textrm{PLF}(\mathbb{S}^{1})= \{f\in
%\textrm{PL}(\mathbb{S}^{1}): B(f)\ \textrm{is finite}\}$.
%
%\textbf{****} \textbullet \  \ $\textrm{PLF}^+(\mathbb{S}^{1})$ the
%group of orientation-preserving elements of
%$\textrm{PLF}(\mathbb{S}^{1})$,

%\textbf{****} \textbullet \  \ $\textrm{PLF}^-(\mathbb{S}^{1})$ the
%set of orientation-reversing elements of
%$\textrm{PLF}(\mathbb{S}^{1})$.
%\medskip

 If $f\in \textrm{Homeo}^+(\mathbb{S}^1)$, we denote by $\rho(f)$ its
rotation number.
\medskip

\emph{In the sequel we identify $\rho(f)$ to its lift in $[0, 1[$.}
\medskip

It is known (see for instance \cite{H}) that if an element $f \in
\textrm{PL}^+(\mathbb{S}^{1})$ is reversed by $h \in
\textrm{PL}^+(\mathbb{S}^{1})$, then by equality (1.1), $\rho(f)=0$
or $\frac{1}{2}$. %The converse is not true (see Example \ref{exe1}).
\vskip 3 mm

We refer the reader to the book \cite{aN10} and papers ( 
%\cite{Ma2}, 
\cite{G}, \cite{Ma}) for a thorough account on groups of circle
homeomorphisms.

The object of this paper is to characterize reversible elements
(resp. strongly reversible elements) in the groups
$\textrm{PL}^+(\mathbb{S}^1)$ and $\textrm{PL}(\mathbb{S}^1)$. Our
main results are the following.
\medskip

 \begin{thm}[Reversibility in $\mathrm{PL}^+(\mathbb{S}^{1})$] \label{t11}
Let $f\in \mathrm{PL}^+(\mathbb{S}^1)$. Then $f$ is reversible in
$\mathrm{PL}^+(\mathbb{S}^1)$ if and only if one of the following
holds:

\begin{enumerate}
  \item[(1)] $\rho(f) = 0$, and $f$ is strongly reversible in $\mathrm{PL}^+(\mathbb{S}^1)$.
  \item[(2)] $\rho(f) = \frac{1}{2}$, and $f$ is strongly reversible in $\mathrm{PL}(\mathbb{S}^{1})$ by an element of 
  $\mathrm{PL}^-(\mathbb{S}^{1})$.
  \end{enumerate}
\end{thm}
\medskip

\begin{rem} \label{r1} If instead of the group $\mathrm{PL}^+(\mathbb{S}^{1})$ we take the group
$\mathrm{Homeo}^+(\mathbb{S}^1)$, the theorem \ref{t11} is false: In
\cite{GOS}, Gill et al. gave an example of a homeomorphism $f\in
\mathrm{Homeo}^+(\mathbb{S}^{1})$ with rotation number $\rho(f) =
\frac{1}{2}$ that is reversible in $\mathrm{Homeo}^+(\mathbb{S}^1)$
but not strongly reversible in $\mathrm{Homeo}(\mathbb{S}^1$).
 \end{rem}
\smallskip

\begin{thm}[Reversibility in
$\mathrm{PL}(\mathbb{S}^{1})$]\label{t2} In
$\mathrm{PL}(\mathbb{S}^{1})$ reversibility and strong reversibility
are equivalent.
\end{thm}
\medskip

We denote by $n_f$ the smallest positive integer $n$ such that $f^n$
has fixed points and by $\Delta_f$ the signature of $f$ (see
definition in Section 2).
\medskip

\begin{thm}\label{t52}
\begin{itemize}
  \item[(1)] Let $f \in  \mathrm{PL}^+(\mathbb{S}^1)$. Then $f$ is strongly reversible in $\mathrm{PL}^+(\mathbb{S}^1)$ if and
    only if one of the following holds.
  \item[(i)]  $f^2 = \mathrm{id}$.
  \item[(ii)] $\mathrm{Fix}(f) \neq \emptyset$ and there exists $h\in \mathrm{PL}^+(\mathbb{S}^1)$ such that
  $\rho(h) = \frac{1}{2}$ and $\Delta_f = -\Delta_f \circ h$.\\
  \item[(2)] Let $f\in \mathrm{PL}^+(\mathbb{S}^1)$. 

\subitem(i) If $\rho(f)\in \Bbb Q$ then $f$ is strongly reversible 
in $\mathrm{PL}(\mathbb{S}^1)$ by an element of $\mathrm{PL}^-(\mathbb{S}^1)$ if and only if
  there exists an involution $h \in \mathrm{PL}^-(\mathbb{S}^1)$ such that $\Delta_{f^{n_f}} = \Delta_{f^{n_f}} \circ h$.
 \subitem(ii) If $\rho(f)\in \mathbb{R}\backslash \Bbb Q$ then $f$ is strongly reversible in $\mathrm{PL}(\mathbb{S}^1)$ by 
 an element of $\mathrm{PL}^-(\mathbb{S}^1)$ if and only if $f$ is conjugate to the rotation $r_{\rho(f)}$ through a
homeomorphism $h$ such that $hrsh^{-1}\in \mathrm{PL}^-(\mathbb{S}^1)$, for some rotation $r$ of $\mathbb{S}^1$.
\end{itemize}
\end{thm}
   \medskip

The next theorem is about composition of reversible (resp.
involution) maps.
\medskip

  \begin{thm}\label{t91}
  We have
  \begin{itemize}
\item[(i)] $\mathrm{PL}^+(\mathbb{S}^{1}) = R_2(\mathrm{PL}^+(\mathbb{S}^{1})) =
  I_3(\mathrm{PL}^+(\mathbb{S}^{1})\neq I_2(\mathrm{PL}^+(\mathbb{S}^{1}))$ and \\
   $R_1(\mathrm{PL}^+(\mathbb{S}^{1}))\subset\neq I_2(\mathrm{PL}(\mathbb{S}^{1}))$.\\
 \item[(ii)] $\mathrm{PL}(\mathbb{S}^1) = R_2(\mathrm{PL}(\mathbb{S}^{1})) = I_3(\mathrm{PL}(\mathbb{S}^{1}))
  \neq I_2(\mathrm{PL}(\mathbb{S}^{1}))$ \\ and  $R_1(\mathrm{PL}(\mathbb{S}^{1})) = I_2(\mathrm{PL}(\mathbb{S}^{1}))
  \neq I_1(\mathrm{PL}(\mathbb{S}^{1}))$.
  \end{itemize}
  \end{thm}
\medskip

The structure of the paper is as follows. In Section 2 we give some
notations and preliminaries results that are needed for the rest of the paper. In
Section 3 we study reversibility in $\textrm{PL}^+(\mathbb{S}^{1})$
of elements $f$ of $\textrm{PL}^+(\mathbb{S}^{1})$ by proving
Theorem \ref{t11}. In Section 4, we study reversibility in
$\textrm{PL}(\mathbb{S}^{1})$ by proving Theorem \ref{t2}. Section 5
is devoted to the characterisation of strong reversibility in
$\textrm{PL}(\mathbb{S}^{1})$ of elements of
$\textrm{PL}^{+}(\mathbb{S}^{1})$. Finally, Section 6 is devoted to
the proof of Theorem \ref{t91}.
\bigskip

%%%%%%%%%%%%%%%%%%%%%%%%%%%%%%%%%%%%%%%%%%%%%%%%%%%%%%%%%%%%%%%%%%%%%%%%%%%%%%%%%%%%%%%%%%%%%%%%%%%%%%%%%%%%%%%%%%%%%%%%%%%%%%%%%%%%%%%%%%%%%%%%%%%%%%%%%%%%%%%%%%%%%%%%%%%%%%%%%%%%%%%%%%%%%%%%%%%%%%%%%%%%%%%%%%%%%%%%%%%%%%%%%%%%%%%%%%%%%%%%%%%%%%%%%%%%%%%%%%%%%%%%%%%%%%%%%%%%%%%%%%%%%%%%%%%%%%%%%%%%%%%%%%%%%%%%%%%%%%%%%%%%%%

\section{\bf Notations and some results}
\medskip

Denote by

\textbullet \ $T$ the translation of $\mathbb{R}$ defined by $T(x) =
x+1$, and for each $a\in \mathbb{R}$ let $T_a$ be the translation
defined by $T_a(x) = x+a$. So $T = T_{1}$.

\textbullet \ \ $r_\pi: z\mapsto -z$  the rotation of
$\mathbb{S}^{1}$ by $\pi$.

\textbullet \ For points $x$ and $y$ in $\mathbb{S}^1$, we denote by
$(x,y)$ the open anticlockwise interval from $x$ to $y$, and by
$[x,y]$ the closure of $(x,y)$. We say that $x<y$ in a proper open
interval $I$ in $\mathbb{S}^1$ if $(x,y)\subset I$.

 \textbullet \ For $f\in \textrm{Homeo}(\mathbb{S}^{1})$ we denote
by  deg$(f)= \begin{cases}
  1, & \textrm{if } f\in \textrm{Homeo}^{+}(\mathbb{S}^{1})\\
 -1, & \textrm{if } f\in \textrm{Homeo}^{-}(\mathbb{S}^{1})
\end{cases}$ the degree of $f$.

\textbullet \ For $f\in \textrm{Homeo}^{+}(\mathbb{S}^{1})$ which
has a fixed point, then each point $x\in \mathbb{S}^{1}$ is either
in $\textrm{Fix}(f)$ or it lies in an open interval component $I$ in
$\mathbb{S}^{1}\backslash \textrm{Fix}(f)$. The \textit{signature} of $f$ (see \cite{GOS}) is
a map $\Delta_f: \mathbb{S}^{1}\longrightarrow \{-1,0,1\}$ given by
$$\Delta_f(x) = \left\{
                           \begin{array}{ll}
                             1, & \textrm{if } f(x)> x \\
                             0, & \textrm{if } f(x)=x \\
                             -1, & \textrm{if } f(x)<x.
                           \end{array}
                         \right.$$
\vskip 2 mm

We have the following lemma.
\bigskip

\begin{lem} \cite{GOS}\label{l100} Let $f\in\mathrm{Homeo}^{+}(\mathbb{S}^{1})$ with a fixed point
and $h \in\mathrm{Homeo}(\mathbb{S}^{1})$. Then\\
$(i)$ $\Delta_{hfh^{-1}} = \mathrm{deg}(h)(\Delta _{f} \circ h^{-1}).$\\
$(ii)$ $\Delta_{ f^{-1}} = -\Delta _{f}.$
\end{lem}
\medskip
\medskip

 \begin{lem}\label{l28}
 Let $f, \ h\in \mathrm{Homeo}(\mathbb{S}^1)$ be such that $f^{-1} = hfh^{-1}$ and let $n\in \mathbb{Z}$.
Then we have $f^{(-1)^{n}} = h^n f h^{-n}$.
\end{lem}

\begin{proof} The proof is down by induction, which is straighforward.
  \end{proof}
  \medskip
  
The following lemma shows that any reversible element of
$\textrm{PL}^+(\mathbb{R})$ must have a fixed point.

\begin{lem}\label{l22}
Let $f, \ h\in \mathrm{PL}^+(\mathbb{R})$ such that
$hfh^{-1} = f^{-1}$. Assume that $f$ is not the identity.
Then $\mathrm{Fix}(f) \neq \emptyset$ and $\mathrm{Fix}(h)
=\emptyset$.
\end{lem}
\medskip

 \begin{proof}
 In fact we show that $f$ has a fixed point in any subinterval $I$ of $\mathbb{R}$ such that $f(I) = I= h(I)$: Otherwise, either $f(x) > x$ for all $x \in I$ or $f(x) < x$ for all $x \in I$. 
 Assume that $f(x) > x$ for all $x\in I$. Then for all $x\in I$, $f^{-1}(x) = hf(h^{-1}(x)) > h(h^{-1}(x))= x$
 since $f$ and $h$ are increasing, which means that $f(x) < x$; a contradiction. Thus $\textrm{Fix}(f_{|I})\neq\emptyset$ and 
in particular, we have $\textrm{Fix}(f)\neq\emptyset$.  As $f\neq \textrm{id}$, so it is not an involution and hence 
$h\neq \textrm{id}$. We will show that $\textrm{Fix}(h)
=\emptyset$. Otherwise,  suppose that $\mathrm{Fix}(h)
\neq\emptyset$. We will prove that in this case $f = \textrm{id}$, this leads to a contradiction. 

Let $I = (a,b)$ be a connected component of $\mathbb{R} \backslash \textrm{Fix}(h)$. 
Then, either $a$ or $b$ is a real number. Let us assume that $a \in \mathbb{R}$. 
By ([2], Lemma 2.4) $f$ fixes  each point of $\textrm{Fix}(h)$. Then $f(I) = I = h(I)$.
Now, by the first step, $f$ has a fixed point $s \in I$. 
So for each integer $n \in \mathbb{Z}$, $h^n(s)\in I$ and by Lemma \ref{l28}, $f((h^n(s)) = h^n(s)$. 
Since $\textrm{Fix}(h_{|I}) = \emptyset$, we can assume, by swapping $h$ and $h^{-1}$ if necessary, that
$h(s) < s < h^{-1}(s)$. Then the points $h^n(s)\in [a, s]$ for $n\in \mathbb{N}$, and accumulate at $a$. 
So, $f$ has infinitely many fixed points in the interval $[a, s]$. 
Since $f \in \textrm{PL}^+(\mathbb{R})$, there is an integer $N$ such that $f_{|[a, h^N(s)]} 
= \textrm{id}$. Thus for any $y\in [a, s]$, one has $f(y)= fh^{-N}(x)$ where $x= h^{N}(y)\in [a,h^{N}(s)]$. By Lemma \ref{l28},
$f(y)= 
h^{-N}f^{(-1)^{N}}(x)=h^{-N}(x)=y$. Therefore $f_{|[a, s]} 
= \textrm{id}$. Similarily, by considering the points $h^{-n}(s) \in [s, b]$, $n\in \mathbb{N}$, we get as above $f_{|[s, b]} 
= \textrm{id}$.  Therefore $f = \textrm{id}$ on $\mathbb{R} \backslash \textrm{Fix}(h)$. As
$f = \textrm{id}$ on $\textrm{Fix}(h)$, thus $f = \textrm{id}$ on $\mathbb{R}$.
\end{proof}
  \medskip

\begin{prop}\label{p24}
\begin{itemize}
  \item[(1)] If $\tau$ is an involution in $\mathrm{PL}^+(\mathbb{S}^{1})$ which is
not the identity then $\tau$ is conjugate in
$\mathrm{PL}^+(\mathbb{S}^{1})$ to the rotation $r_{\pi}$.
\item[(2)] If $\tau$ is an involution in $\mathrm{PL}^-(\mathbb{S}^1)$, then it is conjugate in $\mathrm{PL}(\mathbb{S}^1)$ to 
the reflection $s$.
\end{itemize}
\end{prop}
\medskip

\begin{proof} Let $a$ be the point on $\mathbb{S}^1$ with coordinates $(1, 0)$ and let $b$ be the point with
  coordinates $(-1, 0)$. (1) Let $x \in S^1$. Since $\tau \in
\textrm{PL}^+(\mathbb{S}^{1})$ and $\tau^2 = \textrm{id}$, we have
$\tau([x, \tau(x)]) = [\tau(x), x]$. Let $v : [a, b] \longrightarrow
[x, \tau(x)]$ be a piecewise linear homeomorphism, and let $\psi$ be
the map of $S^1$ defined by $$ \psi(x) = \left\{
                                                             \begin{array}{ll}
                                                               v(x), & \hbox{if} \ x \in [a, b]; \\
                                                               \cr \cr \tau v r_\pi(x), & \hbox{if } \ x \in [b, a].
                                                             \end{array}
                                                           \right.$$ \noindent Then $\psi$ is a well defined piecewise
linear homeomorphism of $S^1$, and it satisfies $\psi^{-1} \tau \psi
= r_\pi$. We conclude that $\tau = \psi r_\pi \psi^{-1}$ is
conjugate in $\textrm{PL}(\mathbb{S}^{1})$ to $r_\pi$.
\smallskip

  \noindent (2)  Let $\{c, d\} = \textrm{Fix}(\tau)$. We have $\tau([c, d]) = [d, c]$.
  Let $u: [a, b] \longrightarrow [c, d]$ be a piecewise
linear homeomorphism and let $\varphi$ be the map of $\mathbb{S}^1$
defined by

$$ \varphi(x) = \begin{cases}
                    u(x), & \textrm{ if }\ x \in [a, b] \\
                    \cr \cr \tau u s(x), & \textrm{ if } \ x \in [b, a]
                  \end{cases}$$ 
                  
\noindent Then $\varphi\in \textrm{PL}(\mathbb{S}^1)$. If $x \in [a, b]$ then 
$\varphi (s(x)) = \tau us(s(x)) = \tau u(x)= \tau \varphi(x)$. If $x \in [b, a]$ then
$\varphi (s(x)) = u(s(x))= \tau \varphi(x)$. Therefore $\tau = \varphi
s\varphi^{-1}$.\end{proof}
 \vskip 3 mm

\vskip 3 mm

 \begin{thm}[\cite{Her}]\label{t26} Let $h\in \mathrm{PL}^+(\mathbb{S}^{1})$ with rotation number $\rho(h)$ irrational.
 Then $h$ is conjugate in $\mathrm{Homeo}^+(\mathbb{S}^{1})$ to the rotation $r_{\rho(h)}$.
\end{thm}
\vskip 3 mm

 \begin{lem}\label{l27}
 Let $f \in \mathrm{Homeo}^+(\mathbb{S}^{1})$ and $g \in \mathrm{Homeo}(\mathbb{S}^{1})$ such that $fg = gf$ and
 $\mathrm{Fix}(f) \neq  \emptyset \neq \mathrm{Fix}(g)$.
 Then $\mathrm{Fix}(f) \cap \mathrm{Fix}(g)\neq \emptyset$.
\end{lem}

\begin{proof} Let $x \in \textrm{Fix}(g)$ and let $\omega_{f}(x)$ be the $\omega$-limit set of $x$ under $f$. It is well known that
$\omega_{f}(x)$ is a periodic orbit (see \cite{H}). Since
$\textrm{Fix}(f) \neq \emptyset$, $\rho(f) = 0$ and any periodic
orbit under $f$ is a fixed point of $f$. It follows that $w_{f}(x) =
\{a\}$; where $a \in \textrm{Fix}(f)$. As $\mathbb{S}^{1}$ is
compact, so $(f^n(x))_{n}$ converges
 to $a$. Now, we have $g(f^n(x)) = f^n(g(x)) = f^n(x)$ since $g$ commutes with $f$. So, $f^n(x)$ converges to $g(a) = a$.
Hence $a\in \textrm{Fix}(f)\cap \textrm{Fix}(g)$. \\
\end{proof}
\vskip 3 mm

  \medskip
  
  %%%%%%%%%%%%%%%%%%%%%%%%%%%%%%%%%%%%%%%%%%%%%%%%%%%%%%%%%%%%%%%%%%%%%%%%%%%%%%%%%%%%%%%%%%%%%%%%%%%%%%%%%%%%%%%%%%%%%%%%%%%%%%%%%%%%%%%%%%%%%%%%%%%%%%%%%%%%%%%%%%%%%%%%%%%%%%%%%%%%%%%%%%%%%%%%%%%%%%%%%%%%%%%%%%%%%%%%%%%%%%%%%%%%%%%%%%%%%%%%%%%%%%%%%%%%%%%%%%%%%%%%%%%%%%%%%%%%%%%%%%%%%%%%%%%%%%%%%%%%%%%%%%%%%%%%%%%%%%%%%%%%%%
  
\section{\bf Reversibility in $\mathrm{PL}^+(\mathbb{S}^{1})$}

\subsection{\bf Reversibility in $\textrm{PL}^{+}(\mathbb{S}^1)$ of
elements $f\in \textrm{PL}^+(\mathbb{S}^{1})$ with $\rho(f)=0$.} The
aim of this subsection is to prove the following proposition.
\medskip

\begin{prop}\label{p21}
 Let $f\in \mathrm{PL}^+(\mathbb{S}^{1})$ such that $\rho(f) = 0$. Then $f$ has a lift $\tilde{f}$ in $\mathrm{PL}^+(\mathbb{R})$,
which is strongly reversible in $\mathrm{PL}(\mathbb{R})$ by an element of
$\mathrm{PL}^-(\mathbb{R})$.
 \end{prop}
\medskip

\begin{proof} Since $\rho(f) = 0$, $f$ has a fixed point $x_0 = e^{it} \in \mathbb{S}^1$. We can assume that $1 \in \textrm{Fix}(f)$
(by taking $r_{-t} f r_t$ instead of $f$, where $r_t(z)= e^{it}z$ is the rotation by angle $t$). Let $\tilde{f}$ be the lift for $f$
such that $\tilde{f}(0) = 0$. Then for all $n \in \mathbb{Z}$, $(T\tilde{f})^n(0) = n$. Let $\alpha_0 : [0, 1] \longrightarrow [0,1]$
be an orientation preserving piecewise linear homeomorphism ($\alpha_0 \in
 \textrm{PL}^+([0,1])$) and for each $n \in \mathbb{Z}$, let $\alpha_{n}: [n, n+1] \rightarrow [n, n+1]$ be the homeomorphism defined as:
$\alpha_n = T^n \alpha_0 (T\tilde{f})^{-n}$. Define $\alpha:
\mathbb{R}\longrightarrow \mathbb{R}$ as $\alpha_{\mid \ [n, n+1]} =
\alpha_n$, for all $n\in \mathbb{Z}$. Then $\alpha \in
\textrm{PL}^+(\mathbb{R})$ and  $T\tilde{f} = \alpha^{-1} T \alpha$.
Moreover, $(T\tilde{f})$ is a lift of $f$. On the other hand, the
translation $T$ satisfies $T^{-1} = iTi$; where $i$ is the
involution of $\mathbb{R}$ defined by $x \longmapsto 1-x$. Then
$(T\tilde{f})^{-1} = \tau (T\tilde{f}) \tau$; where $\tau =
\alpha^{-1}i \alpha$ is an involution in
$\textrm{PL}^-(\mathbb{R})$. The proof is complete.
\end{proof}
\vskip 3 mm

 \begin{prop}\label{p23}
  Let $f\in \mathrm{PL}^+(\mathbb{S}^{1})$ such that $\rho(f) = 0$. If f is reversed by $h \in \mathrm{PL}^+(\mathbb{S}^{1})$, then there exists a lift
$\tilde{f} \in \mathrm{PL}^+(\mathbb{R})$ of $f$ which is reversed in $\mathrm{PL}^+(\mathbb{R})$ by any lift
$\tilde{h}$ of $h$.
  \end{prop}
  \medskip

\begin{proof}
Since $\rho(f) = 0$, $\textrm{Fix}(f) \neq \emptyset$ and there is a
lift $\tilde{f} \in \textrm{PL}^+(\mathbb{R})$ of $f$ such that
$\textrm{Fix}(\tilde{f}) \neq\emptyset$. We have \begin{equation}
f^{-1} = h f h^{-1} \end{equation} \noindent Let $\tilde{h} \in
\textrm{PL}^+(\mathbb{R})$ be a lift of $h$. Then
$\tilde{h}\tilde{f}\tilde{h}^{-1}$ is a lift of $h f h^{-1}$. From
equality (3.1), it follows that $\tilde{f}^{-1} =
\tilde{h}\tilde{f}\tilde{h}^{-1} T_{-q}$ for some integer $q \in
\mathbb{Z}$. As $T \tilde{f} = \tilde{f} T$ and $T \tilde{h} =
\tilde{h} T$,  then $\tilde{f}^{-1} = T_{-q}
\tilde{h}\tilde{f}\tilde{h}^{-1} =  \tilde{h} (T_{-q} \tilde{f})
\tilde{h}^{-1}$. Therefore, since $\tilde{f}$ has a fixed point,
$(T_{-q} \tilde{f})$ has also a fixed point $a\in \mathbb{R}$. Then
$(T_{-q} \tilde{f})(a) = a$, equivalently to $\tilde{f}(a) =
T_q(a)$.
  So, $$\tilde{f}^{n}(a) = {(T_q)}^n(a) = a+nq, \ \ \ \ \ \ \ \forall n \in \mathbb{Z}.$$
\noindent Now, we show that $q = 0$. Suppose that $q \neq 0$ (say $q
> 0$). Then we have $\mathbb{R} = \underset{n \in \mathbb{Z}}\cup\ [a+nq,
a+(n+1)q]$. Let $\alpha_0: [a, a+q] \longrightarrow [0,1]$ be a
piecewise linear orientation preserving homeomorphism such that
$\alpha_0(a) = 0$ and $\alpha_0(a+q) = 1$. For $n \in \mathbb{Z}$,
let $\alpha_n: [a+nq, a+(n+1)q] \rightarrow [n, n+1]$ be a
homeomorphism defined as $\alpha_n= T^n \alpha_0 \tilde{f}^{-n}$.
Define $\alpha: \mathbb{R} \longrightarrow \mathbb{R}$ as
$\alpha_{\mid  [a+nq, a+(n+1)q]} = \alpha_n$ for all $n \in \mathbb{Z}$.
Then  $\alpha\in \textrm{PL}^+(\mathbb{R})$ and we see that
$\tilde{f} = \alpha^{-1}T \alpha$, which is impossible since
$\textrm{Fix}(\tilde{f}) \neq\emptyset$. We conclude that $q=0$ and
so $\tilde{f}^{-1} = \tilde{h}\tilde{f}\tilde{h}^{-1}$.
\end{proof}
\vskip 3 mm

 \begin{prop}\label{t24}
 Let $f\in \mathrm{PL}^+(\mathbb{S}^{1})$ which is not the identity and such that $\rho(f)=0$. Then $f$ is reversible in $\mathrm{PL}^+(\mathbb{S}^{1})$ if and only
 if some lift $\tilde{f}$ of $f$ is reversed by a homeomorphism $\tilde{h}\in
 \mathrm{PL}^+(\mathbb{R})$ satisfying $\tilde{h}T = T\tilde{h}$.
 \end{prop}

\begin{proof} The if part follows from Proposition \ref{p23}. The only if part:  let $f \in \textrm{PL}^+(\mathbb{S}^{1})$ with $\rho(f) = 0$
and let $\tilde{f}$ be a lift for $f$ such that $\tilde{f}^{-1} =
\tilde{h}\tilde{f}\tilde{h}^{-1}$; where $\tilde{h}\in
\textrm{PL}^+(\mathbb{R})$ satisfying $\tilde{h}T = T\tilde{h}$.
Then $\tilde{h}$ is a lift for the homeomorphism $h:
\mathbb{S}^{1}\longrightarrow \mathbb{S}^1$ defined by \ $h(e^{i2
\pi x}) = e^{i2 \pi\tilde{h}(x)}, \ \ \forall \ x \in \mathbb{R}$.
Then $h \in \textrm{PL}^+(\mathbb{S}^1)$ and for all $x \in
\mathbb{R}$, we have $f^{-1} = hfh^{-1}$.
   \end{proof}
 \vskip 3 mm

 \begin{prop}\label{p25} Let $f\in \mathrm{PL}^+(\mathbb{S}^{1})$ such that $\rho(f)=0$.
 If $f$ is a reversible element in $\mathrm{PL}^+(\mathbb{S}^{1})$ then some lift $\tilde{f}$ of $f$ is
conjugate to a homeomorphism $\tilde{g}\in \mathrm{PL}^+(\mathbb{R})$ which is the lift of a homeomorphism $g\in \mathrm{PL}^+(\mathbb{S}^{1})$ that is reversible by the rotation $r_\pi$.
 \end{prop}

 \begin{proof} From Proposition \ref{p23}, there exists $\tilde{h} \in \textrm{PL}^+(\mathbb{R})$ such that
$\tilde{f}^{-1} = \tilde{h}\tilde{f}\tilde{h}^{-1}$. Then by Lemma
\ref{l22}, $\textrm{Fix}(\tilde{h}) =\emptyset$. So $\tilde{h}$ is
conjugate in $\textrm{PL}^+(\mathbb{R})$ to either the translation
$T: x \longmapsto x+1$ or $T_{-1}: x \longmapsto x-1$, say $T$ for
example (see \cite{GS}). It follows that $\tilde{f}$ is conjugate to
a homeomorphism $\tilde{g} \in
 \textrm{PL}^+(\mathbb{R})$ satisfying $\tilde{g}^{-1} = T\tilde{g}T^{-1}$. Therefore $T^2\tilde{g} = \tilde{g}T^2$ and we can define
a homeomorphism $g \in \textrm{PL}^+(\mathbb{S}^{1})$ by \ $g(e^{i \pi x}) = e^{i \pi \tilde{g}(x)}, \ \forall x \in \mathbb{R}$.
We have $g^{-1}(e^{i \pi x}) = e^{i \pi \tilde{g}^{-1}(x)} = e^{i \pi T\tilde{g}T^{-1}(x)} = r_\pi g r_{-\pi}(e^{i \pi x}), \
\forall \ x \in \mathbb{R}$, which means that $g^{-1} = r_\pi g r_\pi$.
   \end{proof}
 \vskip 3 mm

\begin{proof}[Proof of the part (1) of Theorem \ref{t11}]

 \noindent Let $f$ be a reversible homeomorphism in $\textrm{PL}^+(\mathbb{S}^{1})$ such that $\rho(f) =
 0$. If $f$ is the identity, assertion (1) is clear. So suppose that $f$ is not the identity.
Then there exists $h \in PL^+(\mathbb{S}^{1})$ such that $f^{-1} =
hfh^{-1}$. Let us prove that $\rho(h) \in \mathbb{Q}$. Otherwise,
$h$ is conjugate to an irrational rotation $r$ by Theorem \ref{t26};
there exists $\alpha \in \textrm{Homeo}^+(\mathbb{S}^{1})$ such that
$h = \alpha r\alpha^{-1}$. It follows that $f^{-1} = \alpha r
\alpha^{-1} f \alpha r^{-1} \alpha^{-1}$, equivalently $g^{-1} = r g
r^{-1}$; where $g = \alpha^{-1} f \alpha$. Then by Lemma \ref{l28}
for each integer $n \in \mathbb{Z}$,
  \begin{equation} r^{2n}g = g r^{2n}. \end{equation}
  \noindent  Since $\rho(f) = 0$, $\textrm{Fix}(g) \neq \emptyset$. So let $a \in \textrm{Fix}(g)$. By equality (3.2), for
each $n \in \mathbb{Z}$, $r^{2n}(a) \in \textrm{Fix}(g)$. As the
rotation $r$ is irrational, we have $\mathbb{S}^{1}=
\overline{\{r^{2n}(a): \ n \in \mathbb{Z}\}}$. Therefore
$\mathbb{S}^{1}\subset \textrm{Fix}(g)$ and $g = \textrm{id}$. So,
$f = \textrm{id}$, a contradiction.
  \medskip
Let then $\rho(h) = \frac{p}{q}$; where $p$ and $q$ are coprime
positive integers. Then $h^q$ has a fixed point. Let us prove that
$q$ is even.

Otherwise, by Lemma \ref{l28}, $f^{-1} = h^qfh^{-q}$. Then as in the
proof of Proposition \ref{p23}, there is a lift $\tilde{h^{q}}$ of
$h^q$ such that $\textrm{Fix}(\tilde{h^{q}})\neq \emptyset$ and a
lift $\tilde{f}$ for $f$ such that $\tilde{f}^{-1} =
\tilde{h^{q}}\tilde{f}\tilde{h^{q}}^{-1}$; this contradicts Lemma
\ref{l22}. We conclude that $q$ cannot be odd.

Now, one can take $q = 2i$ for some integer $i$. Since $\rho(h^{2i})
= p$, $\textrm{Fix}(h^{2i}) \neq\emptyset$. Since $f h^{2i}=
h^{2i}f$ and $\textrm{Fix}(f) \neq\emptyset$, so
 by Lemma \ref{l27}, there exists $a \in \mathbb{S}^{1}$ such that $f(a) = a = h^{2i}(a)$. It follows that $f(h^j(a)) = h^j(a)$
for each integer $j$ since $f^{-1} = hfh^{-1}$. Let $\mu: \mathbb{S}^{1}\longrightarrow  \mathbb{S}^1$ be the homeomorphism defined by
 $$\mu(z) = \begin{cases}
h^{i}(z), & \textrm{if }  z \in [a, h^i(a)]\\
 h^{-i}(z), & \textrm{if } z \in [h^i(a), a]
 \end{cases}$$
\vskip 3 mm

\noindent It is clear that $\mu$ is an involution in $\textrm{PL}^+(\mathbb{S}^{1})$. We have
 \vskip 2 mm

\begin{equation}\mu f \mu(z) = \begin{cases}
h^{-i}f h^i(z), & \textrm{if} \ z \in [a, h^i(a)]\\
  h^i f h^{-i}(z), & \textrm{if} \ z \in [h^i(a), a]
  \end{cases}
\end{equation}

\noindent We prove that the integer $i$ is odd: Otherwise, suppose
that $i$ is even. Then $f h^i = h^i f$, and by (3.3) we obtain
that $\mu f \mu = f$. Moreover, as $f^{-1} = h f h^{-1}$, then we
have $f^{-1} = h \mu f \mu h^{-1}$.
 Let $a = e^{i2 \pi y}$; where $y\in [0, 1[$, and let $\tilde{h}$ be a lift of $h$.
  As $\rho(h) = \frac{p}{2i}$, then $n = 2i$ is the smallest integer
such that $h^n$ has a fixed point. Therefore, $\tilde{h}^i(y)\neq
y$. Assume that $y < \tilde{h}^i(y)$. Since $\mu(a) = h^i(a)$
  and $\tilde{h}^i$ is a lift for $h^i$, there is a lift $\tilde{\mu}$ for $\mu$ such that $\tilde{\mu}(y) = \tilde{h}^i(y)$.
By Proposition \ref{p23}, there exists a lift $\tilde{f}$ of $f$
such that $\tilde{f}(y) = y$ and $\tilde{f}^{-1} = \tilde{h}
\tilde{\mu} \tilde{f} \tilde{\mu} \tilde{h}^{-1} = \tilde{h}
\tilde{f} \tilde{h}^{-1}$. It follows that $\tilde{\mu} \tilde{f}
\tilde{\mu} = \tilde{f}$. Therefore
${\tilde\mu}\tilde{f}\tilde{\mu}(y) =
{\tilde\mu}\tilde{f}(\tilde{h}^i(y)) = \tilde{f}(y) = y$. Since
$\tilde{\mu} \tilde{f}$ is an increasing homeomorphism, the
inequality  $y< \tilde{h}^i(y)$ implies that $\tilde{\mu}
\tilde{f}(y) \leq y$; equivalently $\tilde{h}^i(y) \leq y$, which is
a contradiction since $y < \tilde{h}^i(y)$. Now the equality $f^{-1}
= hfh^{-1}$ implies that $f^{-1} = h^ifh^{-i}$, equivalently $f^{-1}
= h^{-i}fh^i$. From (3.3), we deduce that $\mu f \mu = f^{-1}$.
Therefore $f$ is strongly reversible by $\mu$ in
$\textrm{PL}^+(\mathbb{S}^{1})$.
  This completes our proof.
  \end{proof}
\vskip 3 mm

\begin{lem}\label{l299} 
  Let $I = (a, b)$ be an open interval in $\mathbb{R}$ or in $\mathbb{S}^{1}$. Then the following statements hold.
\begin{itemize}
\item[(1)] Let $f\in \mathrm{PL}^-(I)$. Then $f$ is reversible in
$\mathrm{PL}(I)$ if and only if $f$ is an involution.
\item[(2)] Let $f\in \mathrm{PL}^+(I)$. Then $f$ is reversible in $\mathrm{PL}(I)$ by an element of
$\mathrm{PL}^-(I)$ if and only if $f$ is strongly
reversible  in $\mathrm{PL}(I)$ by an element of $\mathrm{PL}^-(I)$.
 \end{itemize}
  \end{lem}

\begin{proof} \textbullet \ First, assume that $I$ is an open interval in the real line $\mathbb{R}$.

\noindent Assertion (1). If $f$ is an involution then it is reversible in $\textrm{PL}(I)$ by the identity map. 
Conversely, let $f\in \textrm{PL}^-(I)$, $h \in \textrm{PL}(I)$ such that $hfh^{-1} = f^{-1}$. By replacing $h$ with $hf$ if necessary, 
we can assume that $h \in \textrm{PL}^+(I)$. One can extend $f$ and $h$ on $\mathbb{R}$ as follows: 
$$\hat{f}(x) = \begin{cases}
f(x), & \textrm{if }  x \in (a, b)\\
 a+b-x, & \textrm{if } x \in \mathbb{R} \setminus (a, b)
 \end{cases} \ \ ; \ \ \hat{h}(x) = \begin{cases}
h(x), & \textrm{if }  x \in (a, b)\\
 x, & \textrm{if } x \in \mathbb{R} \setminus (a, b)
 \end{cases}$$
\vskip 3 mm

\noindent Then, clearly $\hat{f} \in \textrm{PL}^-(\mathbb{R})$,  $\hat{h} \in \textrm{PL}^+(\mathbb{R})$ and 
$\hat{h}\hat{f}\hat{h}^{-1} = \hat{f}^{-1}$. So, by ([2], Proposition 2.5), $\hat{f}^2 = \textrm{id}$ on $\mathbb{R}$. 
It follows that $f^2 = \textrm{id}$ on $I$.
\vskip 3 mm

\noindent Assertion (2). Let $f\in \textrm{PL}^+(I)$ and $h\in \textrm{PL}^-(I)$ such that $hfh^{-1} = f^{-1}$. 
Let $\{p\} = \hbox{Fix}(h)$, $p\in I$. We define the involution $\tau \in \textrm{PL}^-(I)$ as follows:
$$ \tau(x) = \begin{cases}
h^{-1}(x), & \textrm{if }  x \geq p\\
 h(x), & \textrm{if } x \leq p
 \end{cases}$$
\vskip 3 mm

\noindent If $f(p) = p$, then clearly $\tau f\tau = f^{-1}$. If  $f(p)\neq p$, let $(c, d)$ be the connected component of 
$I\backslash \hbox{Fix}(f)$ containing $p$ (when $\hbox{Fix}(f) = \emptyset$, \ $(c, d)= I$). 
By the fact that $h f h^{-1} = f^{-1}$, we have $h(\hbox{Fix}(f)) = \hbox{Fix}(f)$. Therefore $h((c, d)) = (c, d) = f((c, d))$ 
(since $h(p) = p$). Let $$ \widehat{f}(x) = \begin{cases}
f(x), & \textrm{if }  x\in (c,d)\\
x, & \textrm{if } x\in \mathbb{R}\backslash (c,d)
 \end{cases}$$ and $$ \widehat{h}(x) = \begin{cases}
h(x), & \textrm{if }  x\in (c,d)\\
x, & \textrm{if } x\in \mathbb{R}\backslash (c,d)
 \end{cases}.$$ Then $\widehat{h}$ is one bump function and satisfies $\hat{h}\hat{f}\hat{h}^{-1} = \hat{f}^{-1}$.
 From ([2], Lemma 4.2), $h^2 = \textrm{id}$ on $(c, d)$ and hence $h(x) = h^{-1}(x)$, for each $x\in (c, d)$. We conclude that 
 the equality $\tau f \tau = f^{-1}$ is satisfied.
\vskip 3 mm

\textbullet \ Now, let us assume that $I = (a, b)$ is an open interval in the circle $\mathbb{S}^{1}$. 
Then, there exists an open interval $\widehat{I} = (t_{1}, t_{2})$ in $\mathbb{R}$ such that the map 
$\varphi: (t_{1}, t_{2})\longrightarrow (a, b)$ given by $\varphi(t) = e^{2i\pi t}$, is a homeomorphism. For 
 $f \in \textrm{PL}(I)$, set $g = \varphi^{-1}f\varphi$. 
We have $g(t) = \varphi^{-1}f\varphi (t) = \varphi^{-1}f(e^{i 2 \pi t}) = 
\varphi^{-1}(e^{i 2\pi\tilde{f}(t)}) = \tilde{f}(t)$, for each $t\in (t_{1}, t_{2})$. Hence $g\in \textrm{PL}(\widehat{I})$.
\vskip 3 mm

 Assertion (1). If $f\in \textrm{PL}^-(I)$ then $g\in \textrm{PL}^{-}(\widehat{I})$.
If $f$ is reversed by $h\in \textrm{PL}(I)$ i.e. $hfh^{-1} = f^{-1}$
 then $g$ is reversed by $k= \varphi^{-1} h \varphi\in \textrm{PL}(\widehat{I})$. Hence by above, $g$ is an involution and so is 
$f$. 
\vskip 3 mm

Assertion (2). If $f\in \textrm{PL}^+(I)$ and $h \in \textrm{PL}^-(I)$ such that $f^{-1} = hfh^{-1}$ then
$g$ is reversed in $\textrm{PL}(\widehat{I})$ by an element of $k = \varphi^{-1}h \varphi\in \textrm{PL}^-(\widehat{I})$. By the above, 
$g$ is strongly reversible in $\textrm{PL}(\widehat{I})$ by an element of $\textrm{PL}^{-}(\widehat{I})$. Hence there exists an 
involution $\tau\in \textrm{PL}^-(\widehat{I})$ such that $\tau g \tau = g^{-1}$. So $(\varphi \tau \varphi^{-1})f
(\varphi \tau \varphi^{-1}) = f^{-1}$. As $(\varphi \tau \varphi^{-1}) \in \textrm{PL}^-(I)$ and is an involution, 
so $f$ is strongly reversible in $\textrm{PL}(I)$ by an element of $\textrm{PL}^-(I)$. 
\end{proof}
\medskip

\begin{lem}\label{l29}
 Let $f \in \mathrm{PL}^+(\mathbb{S}^{1})$ such that $\rho(f) = 0$. If $f$ is strongly reversible in $\mathrm{PL}^+(\mathbb{S}^{1})$,
then it is strongly reversible in $\mathrm{PL}(\mathbb{S}^{1})$ by an element of $\mathrm{PL}^-(\mathbb{S}^{1})$.
 \end{lem}

\begin{proof} If $f$ is strongly reversible in $\textrm{PL}^+(\mathbb{S}^{1})$, then by Proposition \ref{p25},
there exist $\alpha \in \textrm{PL}(\mathbb{S}^{1})$ and $g \in
\textrm{PL}^+(\mathbb{S}^{1})$ such that $f = \alpha g \alpha^{-1}$
and $g^{-1} = r_\pi g r_\pi$. By Proposition \ref{p23}, there exists
a lift $\widetilde{g} \in \textrm{PL}^+(\mathbb{R})$ of $g$ such
that
$$\widetilde{g}^{-1} = T_{\frac{1}{2}}\widetilde{g}T_{\frac{1}{2}}; \ \ \ \ \ (\ast)$$
\noindent Recall that $T_{\frac{1}{2}}$ defined by
$T_{\frac{1}{2}}(t) = t + \frac{1}{2}$ for all $t\in \mathbb{R}$, is
a lift of $r_\pi$.

  Let $\widetilde{s}$ be the involution in $\textrm{PL}^-(\mathbb{R})$ defined
by $\widetilde{s}(t) = -t$. Then $\widetilde{s}T_{\frac{1}{2}}$ is
an involution in $\textrm{PL}^-(\mathbb{R})$, and by equality
$(\ast)$, we have $\widetilde{s} {\widetilde{g}}^{-1} =
(\widetilde{s}T_{\frac{1}{2}})(\widetilde{g}\widetilde{s})(\widetilde{s}T_{\frac{1}{2}})$.
Therefore $(\widetilde{g}\widetilde{s})\in
\textrm{PL}^-(\mathbb{R})$ which is strongly reversible in $\textrm{PL}(\mathbb{R})$ by an element of 
$\textrm{PL}^{-}(\mathbb{R})$. So by Lemma \ref{l299},
$(\widetilde{g}\widetilde{s})^2 = \textrm{id}$; equivalently
${\widetilde{g}}^{-1} = \widetilde{s}\widetilde{g}\widetilde{s}$.
The involution $\widetilde{s}$ satisfies $\widetilde{s}(t+1) =
\widetilde{s}(t) - 1$ for all $t\in \mathbb{R}$. Therefore
$\widetilde{s}$ is a lift for the involution $\sigma$ of
$\textrm{PL}^-(\mathbb{S}^{1})$ defined by
$$\sigma(e^{2i\pi t}) = e^{2i\pi \widetilde{s}(t)}; \ \ \ \ \forall \ t \in
\mathbb{R}.$$ \noindent It follows that for each $t \in \mathbb{R}$,
$g^{-1}(e^{2i\pi t}) = e^{2i\pi
\widetilde{s}\widetilde{g}\widetilde{s}(t)} = \sigma g \sigma(e^{2i
\pi t})$. So $g^{-1} = \sigma g \sigma$. We deduce that $f^{-1} =
\tau f\tau$; where $\tau = \alpha\sigma \alpha^{-1}$ is an
involution in $\textrm{PL}^-(\mathbb{S}^1)$.
\end{proof}
  \medskip

  \medskip

\subsection{\bf Proof of the part (2) of Theorem \ref{t11}}
  We need the following lemma.

\begin{lem}\label{l211}
 Let $f\in \mathrm{PL}^+(\mathbb{S}^{1})$ such that $\mathrm{Fix}(f^n) \neq \emptyset$ for some
  integer $n\in \mathbb{N}^*$.
Then the following are equivalent.
\begin{itemize}
  \item[(1)] $f^n$ is strongly reversible in $\mathrm{PL}(\mathbb{S}^{1})$ by an element of $\mathrm{PL}^-(\mathbb{S}^{1})$.
  \item[(2)] $f$ is strongly reversible in $\mathrm{PL}(\mathbb{S}^{1})$ by an element of 
$\mathrm{PL}^-(\mathbb{S}^{1})$.
\end{itemize}
\end{lem}
\medskip

\begin{proof}  If $f$ is strongly reversible in $\mathrm{PL}(\mathbb{S}^{1})$ by an involution $\mu \in \textrm{PL}^-(\mathbb{S}^{1})$ then $f^{-1} = \mu f \mu$,
which implies that $f^{-n} = \mu f^n \mu$. Conversely, assume that, for some integer $n\in \mathbb{N}^*$,
$f^n$ is strongly reversible in $\mathrm{PL}(\mathbb{S}^{1})$ by an involution $\tau \in
\textrm{PL}^-(\mathbb{S}^{1})$, that is $f^{-n} = \tau f^n \tau$. We will show that $f$ is strongly
reversible in $\mathrm{PL}(\mathbb{S}^{1})$ by an element of $\textrm{PL}^-(\mathbb{S}^{1})$. We have
\begin{equation} (f^{n-1} \tau) f^n (\tau f^{1-n} ) = f^{-n},
\end{equation}
\begin{equation} (f \tau) f^n (\tau f^{-1}) = f^{-n}. \end{equation}
The equality (3.4) implies that  $(f^{n-1} \tau)^2 f^n = f^n
(f^{n-1} \tau)^2$. Therefore from Lemma \ref{l27}, there exists $a
\in \textrm{Fix}(f^n) \cap \textrm{Fix}((f^{n-1} \tau)^2)$. Thus
$f^n(a) = a$, $f^n (\tau(a)) = \tau(a)$ and $f^{n-1} (\tau(a)) =
\tau f^{1-n}(a) = \tau f(a)$. In particular, we have $f
\tau\left((\tau(a), f(a))\right) = (\tau(a), f(a))$. Then the
restriction of $f^n/(\tau(a), f(a))$ is an element of
$\textrm{PL}^+((\tau(a), f(a)))$, which is reversed by $(f \tau)_{|
(\tau(a), f(a))} \in \textrm{PL}^-((\tau(a), f(a)))$ by equality
(3.5). Then by Lemma \ref{l299}, $f^n_{| (\tau(a),
f(a))}$ is strongly reversible in $\textrm{PL}((\tau(a), f(a)))$ by an involution $\sigma \in
\textrm{PL}^-((\tau(a), f(a)))$; that is,
$f^{-n}_{| (\tau(a),
f(a))} = \sigma f^n_{| (\tau(a), f(a))} \sigma$. The point $f(a)$ is
either in $(a, \tau(a))$ or in $(\tau(a), a)$. We can assume that
$f(a) \in (\tau(a), a)$. Since $f$ is orientation-preserving, we can
easily see that
$$\mathbb{S}^{1} = \bigcup_{p=1}^n \ [f^p(a), f^p \tau(a)] \cup \bigcup_{p=2}^{n+1}\ [f^p \tau(a), f^{p+1}(a)].$$
\noindent Let $\mu: \mathbb{S}^{1}\longrightarrow \mathbb{S}^{1}$ be
the map of $\mathbb{S}^{1}$ defined by
$$\mu (x) = \begin{cases}
f^{n-p+1}\tau f^{n-p}(x), & \textrm{if }  x \in [f^p(a), f^p\tau(a)], \ \ \forall \ 1 \leq p \leq n; \\
f^{n-p} \sigma f^{n-p}(x), & \textrm{if } x \in (f^p\tau(a), f^{p+1}(a)), \ \ \forall \ 2 \leq p \leq n+1
\end{cases}$$

 \noindent The map $\mu$ is a well defined homeomorphism of $\mathbb{S}^1$. Moreover
 $\mu\in \textrm{PL}^-(\mathbb{S}^{1})$ that satisfies
  $\mu^2 = \textrm{id}$ and $\mu f \mu = f^{-1}$. This completes the proof.
  \end{proof}
\vskip 3 mm

\noindent \textit{Proof of the part (2) of Theorem \ref{t11} } Let
$f\in \textrm{PL}^+(\mathbb{S}^{1})$ be reversible in
$\textrm{PL}^+(\mathbb{S}^{1})$. Then the rotation number $\rho(f)$
is equal to either $0$ or $\frac{1}{2}$. The first case $\rho(f) =
0$ corresponds to the first part of Theorem \ref{t11}. In the second
case, $\rho(f) = \frac{1}{2}$ we have $\rho(f^2) = 0 (\textrm{mod }
1)$. Let us prove that $f$ is strongly reversible in $\textrm{PL}(\mathbb{S}^{1})$ by an element of
$\textrm{PL}^-(\mathbb{S}^{1})$. By hypothesis, there exists a
homeomorphism $h \in \textrm{PL}^+(\mathbb{S}^{1})$ such that
$f^{-1} = hfh^{-1}$ and so $f^{-2} = h f^2 h^{-1}$. Then by the
proof of the part (1) of Theorem \ref{t11}, we know that either
$\rho(h) \in \mathbb{R} \setminus \mathbb{Q}$  or $\rho(h) =
\frac{1}{2i}$ (mod 1); where $i$ is an odd integer. It follows that
in the first case, $f^2 = \textrm{id}$, and in the second case,
$f^2$ is strongly reversible in $\textrm{PL}^+(\mathbb{S}^{1})$ (see
the proof of the part (1) of Theorem \ref{t11}). By Lemma \ref{l29},
$f^2$ is strongly reversible in $\textrm{PL}(\mathbb{S}^{1})$ by an element of
 $\textrm{PL}^-(\mathbb{S}^{1})$. As $\textrm{Fix}f^{2}\neq \emptyset$, we conclude by Lemma \ref{l211} that $f$ is 
 strongly reversible in $\textrm{PL}(\mathbb{S}^{1})$ by an element of $\textrm{PL}^-(\mathbb{S}^{1})$. \ \qed \\
\medskip

%%%%%%%%%%%%%%%%%%%%%%%%%%%%%%%%%%%%%%%%%%%%%%%%%%%%%%%%%%%%%%%%%%%%%%%%%%%%%%%%%%%%%%%%%%%%%%%%%%%%%%%%%%%%%%%%%%%%%%%%%%%%%%%%%%%%%%%%%%%%%%%%%%%%%%%%%%%%%%%%%%%%%%%%%%%%%%%%%%%%%%%%%%%%%%%%%%%%%%%%%%%%%%%%%%%%%%%%%%%%%%%%%%%%%%%%%%%%%%%%%%%%%%%%%%%%%%%%%%%%%%%%%%%%%%%%%%%%%%%%%%%%%%%%%%%%%%%%%%%%%%%%%%%%%%%%%%%%%%%%%%%%%%

\section{\bf Reversibility in $\textrm{PL}(\mathbb{S}^{1})$}
\medskip

\subsection{\bf Reversibility in $\textrm{PL}(\mathbb{S}^{1})$ of
elements of $\textrm{PL}^+(\mathbb{S}^{1})$} The aim of this section
is to prove the following proposition.
\medskip

\begin{prop}\label{p13}
    Let $f \in \mathrm{PL}^+(\mathbb{S}^{1})$. Then $f$ is reversible in $\mathrm{PL}(\mathbb{S}^{1})$ by an element of 
    $\mathrm{PL}^-(\mathbb{S}^{1})$ if and only if it is strongly reversible in $\mathrm{PL}(\mathbb{S}^{1})$ by an element of 
    $\textrm{PL}^-(\mathbb{S}^{1})$.
  \end{prop}
\medskip

  \begin{lem}\label{l31} \begin{itemize}
\item[(1)] Let $I = (a, b)$ be an open interval in $\mathbb{R}$ or in $\mathbb{S}^{1}$. 
Then every fixed point free element $v\in \mathrm{PL}^+(I)$ is strongly 
reversible in $\mathrm{PL}(I)$ by an element of $\mathrm{PL}^-(I)$.

 \item[(2)] Let $f \in \mathrm{PL}^+(\mathbb{S}^{1})$. If $f$ has exactly one
fixed point, then $f$ is strongly reversible in $\mathrm{PL}(\mathbb{S}^{1})$ by an element of $\mathrm{PL}^-(\mathbb{S}^{1})$.
 \end{itemize}
\end{lem}

\begin{proof} (1) \textbullet \ Let $I$ be an open interval in $\mathbb{R}$ and $v\in \textrm{PL}^+(I)$ be a 
fixed point free homeomorphism. A similar construction as in the proof of (\cite{J}, Theorem 1)
prove that there exists an involution $\alpha \in  \textrm{PL}^-(I)$ satisfying  $v^{-1} = \alpha v\alpha$. 

\textbullet \ Assume that $I = (a, b)$ is an open interval in $\mathbb{S}^{1}$. 
Then, there exists an open interval $\widehat{I} = (t_{1}, t_{2})$ in $\mathbb{R}$ such that the map 
$\varphi: (t_{1}, t_{2})\longrightarrow (a, b)$ given by $\varphi(t) = e^{2i\pi t}$, is a homeomorphism. If $v$ is a fixed point free element in 
$\textrm{PL}^+(I)$, then $\varphi^{-1}v\varphi$ is a fixed point free element in $\textrm{PL}^+(\widehat{I})$. Then, by the above,
 there exists an involution $\alpha \in \textrm{PL}^-(\widehat{I})$ satisfying 
$\varphi^{-1} v^{-1} \varphi = \alpha \varphi^{-1}v\varphi \alpha $. 
It follows that $v^{-1} = (\varphi \alpha \varphi^{-1}) v (\varphi \alpha \varphi^{-1})$. As 
$\tau = \varphi \alpha \varphi^{-1} \in \textrm{PL}^-(I)$, we conclude that $v$ is strongly reversible in $ \textrm{PL}(I)$ 
by an element of $\textrm{PL}^-(I)$.
\medskip

(2) Let $\{a\} = \textrm{Fix}(f)$. By  (1), the restriction $f_{| \ \mathbb{S}^{1}\setminus \{a\}}$ 
 is strongly reversible in $\textrm{PL}(\mathbb{S}^{1}\setminus \{a\})$ by an involution $\sigma\in
\textrm{PL}^-(\mathbb{S}^{1}\setminus \{a\})$. Then we extend
$\sigma$ to a map $\hat{\sigma}: \mathbb{S}^{1}\longrightarrow
\mathbb{S}^{1}$
 given by $$\hat{\sigma}(x) = \begin{cases} \sigma(x), & \textrm{if } \ x \in \mathbb{S}^{1}\setminus \{a\},\\
 a, & \textrm{if }\ x = a
\end{cases}$$
\vskip 3 mm

 \noindent We see that $\hat{\sigma}$ is an involution in $\textrm{PL}^-(\mathbb{S}^{1})$ which satisfies
$f^{-1} = \hat{\sigma}f\hat{\sigma}$.
   \end{proof}
 \medskip

 \begin{lem}\label{l32} Let $f\in \textrm{PL}^+(\mathbb{S}^{1})$ such that $\rho(f) = 0$. If $f$ is reversible in
 $\mathrm{PL}(\mathbb{S}^{1})$ by an element of $\mathrm{PL}^-(\mathbb{S}^{1})$ then it is strongly reversible in 
 $\mathrm{PL}(\mathbb{S}^{1})$ by an element of 
 $\textrm{PL}^-(\mathbb{S}^{1})$.
  \end{lem}
  \medskip

 \begin{proof} Assume that there exists $h \in \textrm{PL}^-(\mathbb{S}^{1})$ such that $f^{-1} = h f h^{-1}$.
 Let us show that $f$ is strongly reversible in $\textrm{PL}(\mathbb{S}^{1})$ by an element of $\textrm{PL}^-(\mathbb{S}^{1})$.
   If $f$ has exactly one fixed point, then the conclusion follows from Lemma \ref{l31}. Now, assume that $f$ has more than one fixed point. Since $h \in
\textrm{PL}^-(\mathbb{S}^{1})$, so $h$ has exactly two fixed points
$a$ and $b$ which divides the circle $\mathbb{S}^{1}$ onto two
connected components $A=(a,b)$ and $B=(b,a)$ satisfying $h(A) = B$
and $h(B) = A$. Moreover we have always $\textrm{Fix}(f) \cap A \neq
\emptyset$ and $\textrm{Fix}(f) \cap B \neq \emptyset$. Let $c$ be
the nearest point of $\textrm{Fix}(f) \cap A$ to the point $a$, and
let $d$ be the nearest point of $\textrm{Fix}(f) \cap A$ to the
point $b$. From the equality $f^{-1} = h f h^{-1}$, we see that
$h(c)$ is the nearest point of $\textrm{Fix}(f) \cap B$ to $a$ and
that $h(d)$ is the nearest point of $\textrm{Fix}(f) \cap B$ to $b$.
The restrictions $f_{|(h(c), c)}$ and $f_{|(d, h(d))}$ are fixed
point free piecewise linear homeomorphisms of open arcs of the
circle $\mathbb{S}^{1}$. Then by
 Lemma \ref{l31}, (1), they are reversed respectively in $\textrm{PL}((h(c), c))$ and $\textrm{PL}((d, h(d)))$ by involutions
$\sigma_1 \in \textrm{PL}^-((h(c), c))$ and $\sigma_2 \in
\textrm{PL}^-((d, h(d)))$.  Define $\tau: \mathbb{S}^{1}\longrightarrow   \mathbb{S}^{1}$ as follows
\hskip 53 mm $$\ \  \tau(x)= \begin{cases}
 h(x), & \textrm{if }  x \in [c,d],\\
 h^{-1}(x), & \textrm{if } x \in [h(d), h(c)]\\
 \sigma_1(x), & \textrm{if } x\in (h(c), c)\\
 \sigma_2(x), & \textrm{if } x\in (d, h(d))
\end{cases}$$
\vskip 3 mm

 \noindent We can easily see that $\tau$ is an involution in $\textrm{PL}^-(\mathbb{S}^{1})$ that satisfies $f^{-1} = \tau f \tau$.
  \end{proof}

  \noindent  \emph{\it Proof of Proposition \ref{p13}}. Assume that $f$ is reversible in $\textrm{PL}(\mathbb{S}^{1})$ 
  by an element of 
  $\textrm{PL}^-(\mathbb{S}^{1})$.
  We distinguish two cases:

  \noindent Case 1: $\rho(f) = 0$. Then $f$ is strongly reversible in $\textrm{PL}(\mathbb{S}^{1})$ by an element of $\textrm{PL}^-(\mathbb{S}^{1})$ by Lemma \ref{l32}.

  \noindent Case 2: $\rho(f) \in \mathbb{R} \setminus  \mathbb{Q}$. Then by Theorem \ref{t26}, there exists
  $\alpha \in \textrm{Homeo}(\mathbb{S}^{1})$ such that $f = \alpha r \alpha^{-1}$; where $r$ is the rotation of $\mathbb{S}^{1}$
  by $\rho(f)$. On the other hand, there exists $h \in \textrm{PL}^-(\mathbb{S}^{1})$ such that $f^{-1} = h f h^{-1}$, which implies that $r^{-1} = g r g^{-1}$, where $g = \alpha^{-1}h \alpha$. Then \begin{equation}g^2r = rg^2. \end{equation}
  \noindent Since $g$ is an orientation-reversing element of $\textrm{Homeo}(\mathbb{S}^{1})$,
  $\textrm{Fix}(g) \neq\emptyset$. Let $a \in \textrm{Fix}(g) \subset \textrm{Fix}(g^2)$. The equality
(4.1) implies that for each $n\in \mathbb{Z}$, $g^2 r^n = r^n g^2$.
It follows that $r^n(a) \in \textrm{Fix}(g^2)$, for each $n\in
\mathbb{Z}$ and by the fact that $\mathbb{S}^{1} =
\overline{\{r^n(a): \ n\in \mathbb{Z}\}}$, we obtain that
$\mathbb{S}^{1} = \textrm{Fix}(g^2)$. Thus $g^2 = \textrm{id}$ and
so $h^2 = \textrm{id}$. We conclude that $f$ is strongly reversible in $\textrm{PL}(\mathbb{S}^{1})$ by the involution $h\in \textrm{PL}^-(\mathbb{S}^{1})$.

  \noindent Case 3. $\rho(f) = \frac{p}{q}\in \mathbb{Q}\backslash \{0\}$. In this case, $\rho(f^q) = 0$ and by the case 1,
$f^q$ is strongly reversible in $\textrm{PL}(\mathbb{S}^{1})$ by an element of $\textrm{PL}^-(\mathbb{S}^{1})$. So
by Lemma \ref{l211}, $f$ is strongly reversible in
$\textrm{PL}(\mathbb{S}^{1})$ by an element of $\textrm{PL}^-(\mathbb{S}^{1})$. \qed \vskip 3 mm
\medskip

\subsection{\bf Reversibility in $\textrm{PL}(\mathbb{S}^{1})$ of
elements of $\textrm{PL}^-(\mathbb{S}^{1})$} \ \  In this paragraph
we study reversibility in $\textrm{PL}(\mathbb{S}^{1})$ of elements of
$\textrm{PL}^-(\mathbb{S}^{1})$ by
proving the following proposition.

\begin{prop}\label{p12}
    Let $f \in \mathrm{PL}^-(\mathbb{S}^{1})$. Then the following statements are equivalent.
\begin{enumerate}
  \item $f$ is reversible in $\textrm{PL}(\mathbb{S}^{1})$ by an element of $\mathrm{PL}^+(\mathbb{S}^{1})$.
   \item $f$ is reversible in $\textrm{PL}(\mathbb{S}^{1})$ by an element of $\mathrm{PL}^-(\mathbb{S}^{1})$.
    \item $f$ is strongly reversible in $\textrm{PL}(\mathbb{S}^{1})$ by an element of $\mathrm{PL}^-(\mathbb{S}^{1})$.
     \item $f$ is strongly reversible in $\textrm{PL}(\mathbb{S}^{1})$ by an element of $\mathrm{PL}^+(\mathbb{S}^{1})$.
  \end{enumerate}
  \end{prop}
\medskip

  \begin{lem}\label{l61}
    Let $f \in \mathrm{PL}^-(\mathbb{S}^{1})$. Then the following statements are equivalent.
     \begin{enumerate}
  \item $f$ is reversible by an element $h \in \mathrm{PL}(\mathbb{S}^{1})$ that fixes each of the fixed points of $f$.
 \item $f^2 = \mathrm{id}$.
 \end{enumerate}
  \end{lem}

  \begin{proof}  $(1)\Longrightarrow (2)$:  Since $f$ is orientation-reversing, so it has exactly two fixed points
  $a$ and $b$. Set $I = \mathbb{S}^{1} \backslash \{a\}$, it is an open interval in $\mathbb{S}^{1}$.
  As $f^{-1} = h f h^{-1}$ and $h(a) = a$, the restriction $f_{| I} \in \hbox{PL}^-(I)$ and is reversed by 
$h_{| I} \in \hbox{PL}(I)$. Then by Lemma \ref{l299}, (1),
 $f_{| I}$ is an involution and so is $f$.

   $(2)\Longrightarrow (1)$: is clear.
   \end{proof}
   \vskip 3 mm

  \noindent  \emph{\it Proof of Proposition \ref{p12}.} \noindent $(1) \Longrightarrow (2)$: Let $f\in
\textrm{PL}^-(\mathbb{S}^{1})$, $h \in \textrm{PL}^+(\mathbb{S}^{1})$ such that $f^{-1} = h f
  h^{-1}$. So $f^{-1} = (fh) f(h^{-1}f^{-1})$. Hence $f$ is reversed in $\textrm{PL}(\mathbb{S}^{1})$ by 
  $fh\in \textrm{PL}^-(\mathbb{S}^{1})$.

 \noindent  $(2) \Longrightarrow (3)$: Let $f\in \textrm{PL}^-(\mathbb{S}^{1})$ and
 $h \in \textrm{PL}^-(\mathbb{S}^{1})$ be such that $$f^{-1} = h f h^{-1}. \ \ \ \ \ (\ast)$$
  Since $f$ is orientation-reversing, so it has exactly two fixed points
  $a$ and $b$. We have $h(\textrm{Fix}(f)) = \textrm{Fix}(f)$. So, either $h$ fixes each of $a$ and $b$ or it
  interchanges them.
  In the first case, we have $f^2 = \textrm{id}$ by Lemma \ref{l61}. So
  $f$ is an involution in $\textrm{PL}^-(\mathbb{S}^{1})$ and hence
   it is strongly reversible in $\textrm{PL}(\mathbb{S}^{1})$ by an element of $\textrm{PL}^-(\mathbb{S}^{1})$. In the second case; that is $h(a)
= b$ and $h(b) = a$, we have $h((a,b)) = (a,b)$. So by
  equality $(\ast)$, the restriction $f^2_{|(a,b)}$ is an element of $\textrm{PL}^+((a,b))$ that is reversed by
  $h_{|(a,b)}$. Thus, by Lemma \ref{l299}, (2), $f^2_{|(a,b)}$ is strongly reversible in $\textrm{PL}((a,b))$ 
  by an involution $\tau\in
  \textrm{PL}^-((a,b))$; that is, $$f^{-2}_{|(a,b)} = \tau f^2_{|(a,b)} \tau. \ \ \ \ \ \ \ \ (\ast \ast)$$
  \noindent Let $\mu: \mathbb{S}^{1} \longrightarrow \mathbb{S}^{1}$ be the map defined by  $$\mu(x) = \left\{
                                                                                 \begin{array}{ll}
                                                                                   \tau(x), &  \textrm{if } \ x \in [a, b] \\
                                                                                    f^{-1} \tau f^{-1}(x), & \textrm{if
                                                                                    }
                                                                                    \
                                                                                                              x \in [b, a]
                                                                                     \end{array}
                                                                                       \right.$$
  \noindent Clearly $\mu \in \textrm{PL}^-(\mathbb{S}^{1})$ and $\mu f \mu = f^{-1}$. Moreover,
  by equality $(\ast \ast)$, we have $\mu^2 = \textrm{id}$. This
  implies that $f$ is strongly reversible in $\textrm{PL}(\mathbb{S}^{1})$ by an element of 
  $\textrm{PL}^-(\mathbb{S}^{1})$.

 \noindent  $(3) \Longrightarrow (4)$. Assume that $f^{-1} = \tau f \tau$, where $\tau$ is an involution in
 $\textrm{PL}^-(\mathbb{S}^{1})$. Then $(f \tau)^2 = \textrm{id}$ and so $f^{-1} = (f\tau) f (\tau
f^{-1})$. Hence $f$ is also strongly reversible in $\textrm{PL}(\mathbb{S}^{1})$ by the
  involution $(f\tau)$ in $\textrm{PL}^+(\mathbb{S}^{1})$.

 \noindent $(4) \Longrightarrow (1)$ is clear. \qed
   \bigskip

{\it Proof of Theorem \ref{t2}}. This follows from Theorem
\ref{t11}, Propositions \ref{p13} and \ref{p12}. \qed
\medskip

\begin{rem} \label{r3}  In Proposition
\ref{p12}, we showed that any reversible element in $\textrm{PL}(\mathbb{S}^{1})$ by an element of 
$\textrm{PL}^-(\mathbb{S}^{1})$ must be strongly reversible. This
does not hold for elements of $\textrm{Homeo}^-(\mathbb{S}^1)$ (see
\cite{GOS}).
\end{rem}
\vskip 3 mm

%%%%%%%%%%%%%%%%%%%%%%%%%%%%%%%%%%%%%%%%%%%%%%%%%%%%%%%%%%%%%%%%%%%%%%%%%%%%%%%%%%%%%%%%%%%%%%%%%%%%%%%%%%%%%%%%%%%%%%%%%%%%%%%%%%%%%%%%%%%%%%%%%%%%%%%%%%%%%%%%%%%%%%%%%%%%%%%%%%%%%%%%%%%%%%%%%%%%%%%%%%%%%%%%%%%%%%%%%%%%%%%%%%%%%%%%%%%%%%%%%%%%%%%%%%%%%%%%%%%%%%%%%%%%%%%%%%%%%%%%%%%%%%%%%%%%%%%%%%%%%%%%%%%%%%%%%%%%%%%%%%%%%%

\section{\bf Strong reversibility in $\textrm{PL}(\mathbb{S}^{1})$ of elements of $\textrm{PL}^+(\mathbb{S}^{1})$}
\medskip

\subsection{Strong reversibility of elements of
$\textrm{PL}^+(\mathbb{S}^{1})$}
\medskip

The aim of this subsection is to prove the part (1) of Theorem
\ref{t52}.
\medskip

 \begin{lem}\label{l53}
    Let $f, \ g \in \mathrm{PL}^+(\mathbb{S}^1)$ such that $\mathrm{Fix}(f) \neq \emptyset \neq \mathrm{Fix}(g)$.
    If $\Delta_f = \Delta_g$, then there exists $v \in  \mathrm{PL}^+(\mathbb{S}^1)
    $ such that $g = vfv^{-1}$ and $v = \mathrm{id}$ on $\mathrm{Fix}(f)$.
  \end{lem}

 \begin{proof} Since $\Delta_f = \Delta_g$, we have $\textrm{Fix}(f) = \textrm{Fix}(g)$.
  For each open interval component $(a,b)$ of $S^1 \backslash \textrm{Fix}(f)$,
  there exists an orientation preserving piecewise linear homeomorphism
  $u: (a,b) \longrightarrow \mathbb{R}$. Then $u f u^{-1}$ and $u g u^{-1}$ are two fixed point free elements
  of $\textrm{PL}^+(\mathbb{R})$. Since $\Delta_f = \Delta_g$, $u f u^{-1}$ and $u g u^{-1}$
 are conjugate in $\textrm{PL}^+(\mathbb{R})$ by (\cite{GS}, Proposition 2.6). Let $v_0 \in PL^+((a,b))$ such that $g(x) = v_0 f v_0^{-1}(x)$ for $x \in (a, b)$. Then the map $v$ defined by
$v(x) = v_0(x)$ for $x \in (a, b)$ and $v(x) = x$ for $x \in
\textrm{Fix}(f)$, is the required homeomorphism.
\end{proof}
\medskip

    \begin{lem}\label{l54}
 Let $f \in \mathrm{PL}^+(\mathbb{S}^1)$ such that $\rho(f) = \frac{1}{2}$. Then $f$ is strongly reversible in
 $\mathrm{PL}^+(\mathbb{S}^1)$ if and only if $f^2 = \mathrm{id}$.
 \end{lem}
\medskip

\begin{proof} Lemma \ref{l54} is a particular case of (\cite{GOS}, Theorem 3.3).
\end{proof}
\medskip

   \noindent  \emph{Proof of the part (1) of Theorem \ref{t52}}.
Assume that $f$ is a strongly reversible element of $\textrm{PL}^+(\mathbb{S}^1)$. We know that either $\rho(f) = 0$ or $\rho(f) = {1 \over 2}$.
 If $\rho(f) = {1 \over 2}$, then by Lemma \ref{l53}, $f^2 = \textrm{id}$. If $\rho(f) = 0$, then
$\textrm{Fix}(f)\neq \emptyset$, and since $f$ is strongly reversible in $\textrm{PL}^+(\mathbb{S}^1)$, there exists
 an involution $h \in \textrm{PL}^+(\mathbb{S}^1)$ such that $f^{-1} = h^{-1} f h$. Therefore,
$\rho(h) = \frac{1}{2}$ and by (\cite{GOS}, Lemma 2.1), $\Delta_f =
-\Delta_f \circ h$.

Conversely, assume that $f \in
\textrm{PL}^+(\mathbb{S}^1)$ such that $\textrm{Fix}(f)\neq
\emptyset$ and there exists $h\in \textrm{PL}^+(\mathbb{S}^1)$ with
$\rho(h) = \frac{1}{2}$ satisfying $\Delta_f = -\Delta_f \circ h$. Then $\Delta_{f^{-1}} = \Delta_{h^{-1} f h}$.
By Lemma \ref{l53}, there exists $v \in \textrm{PL}^+(\mathbb{S}^1)$ such that $f^{-1} = v^{-1}h^{-1} f h v$;
which means that $f$ is reversible by $hv\in \textrm{PL}^+(\mathbb{S}^1)$. Since $\textrm{Fix}(f)\neq \emptyset$,
 $\rho(f) = 0$ and by Theorem \ref{t11}, $f$ is strongly reversible in $\textrm{PL}^+(\mathbb{S}^1)$.
If $f^2 = \textrm{id}$, then it is clear that $f$ is strongly reversible in $\textrm{PL}^+(\mathbb{S}^1)$
by the identity map. \qed \\

\subsection{Strong reversibility of elements of
$\textrm{PL}^-(\mathbb{S}^{1})$}
\medskip

The aim of this subsection is to prove the part (2) of Theorem
\ref{t52}.
\medskip

   \begin{lem}\label{l55}
    Let $f \in \mathrm{PL}^+(\mathbb{S}^1)$ with $\rho(f) = 0$. Then $f$ is strongly reversible in $\mathrm{PL}(\mathbb{S}^1)$
by an element of $\mathrm{PL}^-(\mathbb{S}^1)$ if and only if there exists $h \in
\mathrm{PL}^-(\mathbb{S}^1)$ such that $\Delta_f = \Delta_f \circ
h$.
  \end{lem}

  \noindent  \emph{Proof.} Let $f \in \textrm{PL}^+(\mathbb{S}^1)$ such that $\rho(f) = 0$. If $f$
  is strongly reversible in $\mathrm{PL}(\mathbb{S}^1)$ by an element of $\mathrm{PL}^-(\mathbb{S}^1)$ then there exists an involution $h \in \textrm{PL}^-(\mathbb{S}^1)$ such that
  $f^{-1} = h^{-1}fh$. Thus by Lemma \ref{l100}, $\Delta_f = -\textrm{deg}(h) \Delta_f \circ h = \Delta_f \circ h$.
  Conversely, if $\Delta_f = \Delta_f \circ h$ for some element $h \in \textrm{PL}^-(\mathbb{S}^1)$ then
$\Delta_f = -\textrm{deg}(h) \Delta_f \circ h$ and $\Delta_{f^{-1}}
= \Delta_{h^{-1} f h}$. So, by the proof of the part (1) of Theorem
\ref{t52}, $f$ is reversible in $\textrm{PL}(\mathbb{S}^1)$ by an element of $\textrm{PL}^-(\mathbb{S}^1)$. By
Proposition \ref{p13}, $f$ is strongly reversible in $\textrm{PL}(\mathbb{S}^1)$ by an element of 
$\textrm{PL}^-(\mathbb{S}^1)$. \qed
\bigskip

   \noindent  \emph{Proof of the part (2) of Theorem \ref{t52}}. (i): Let $f \in \textrm{PL}^+(\mathbb{S}^1)$ 
   such that $\rho(f) \in \Bbb Q$. If $f$ is strongly reversible in $\textrm{PL}(\mathbb{S}^1)$ by an element of 
   $\textrm{PL}^-(\mathbb{S}^1)$, then there exists an involution $h \in \textrm{PL}^-(\mathbb{S}^1)$ such that
    $f^{-1} = h^{-1}fh$, which implies that $f^{-n_f} = h^{-1} f^{n_f} h$. So by Lemma \ref{l55}, we have
    $\Delta_{f^{n_f}} = \Delta_{f^{n_f}} \circ h$. Conversely, assume that there exists
    $h \in \textrm{PL}^-(\mathbb{S}^1)$ such that $\Delta_{f^{n_f}} = \Delta_{f^{n_f}} \circ h$. Then by Lemma \ref{l55},
     $f^{n_f}$ is strongly reversible in $\textrm{PL}(\mathbb{S}^1)$ by an element of $\textrm{PL}^-(\mathbb{S}^1)$ and by 
     Lemma \ref{l211}, 
     $f$ is strongly reversible in $\textrm{PL}(\mathbb{S}^1)$ by an element of $\textrm{PL}^-(\mathbb{S}^1)$. 
     
     (ii): Let $f \in \textrm{PL}^+(\mathbb{S}^1)$ such that $\rho(f)\in \mathbb{R}\setminus \Bbb Q$. If $f$ is 
     strongly reversible in $\textrm{PL}(\mathbb{S}^1)$ by an element of $\textrm{PL}^-(\mathbb{S}^1)$ then 
   there exists an involution $\tau\in \textrm{PL}^-(\mathbb{S}^1)$ such that
  $f^{-1} = \tau f\tau$.  On the other hand, by Theorem \ref{t26}, there is $h\in \textrm{Homeo}^+(\mathbb{S}^{1})$ 
  such that $f = hr_{\rho(f)}h^{-1}$. Therefore $hr_{\rho(f)}^{-1}h^{-1}= \tau hr_{\rho(f)} h^{-1}\tau$. As $r_{\rho(f)}^{-1} = sr_{\rho(f)}s$, where
  $s: z\mapsto \overline{z}$ is the reflection, so 
  $h^{-1}\tau h sr_{\rho(f)}= r_{\rho(f)}h^{-1}\tau h s$. Hence $h^{-1}\tau h s = r_{t}$ for some $t \in \mathbb{R}$. 
  It follows that 
  $\tau = hr_{t}sh^{-1}\in \textrm{PL}^-(\mathbb{S}^1)$. 
  
  Conversely, if there is $h\in \textrm{Homeo}^+(\mathbb{S}^{1})$ 
  such that $f = hr_{\rho(f)}h^{-1}$ where $h$ satisfies $hrsh^{-1}\in \textrm{PL}^-(\mathbb{S}^1)$, for some rotation $r$ of 
  $\mathbb{S}^1$, then $\tau= hrsh^{-1}$ is an involution in $\textrm{PL}^-(\mathbb{S}^1)$ and satisfies $f^{-1} = \tau f\tau$ 
  (since  $rs = sr^{-1}$). Therefore
  $f$ is strongly reversible in $\textrm{PL}(\mathbb{S}^1)$ by an element of $\textrm{PL}^-(\mathbb{S}^1)$. \hfill \qed \\
\vskip 3 mm

%%%%%%%%%%%%%%%%%%%%%%%%%%%%%%%%%%%%%%%%%%%%%%%%%%%%%%%%%%%%%%%%%%%%%%%%%%%%%%%%%%%%%%%%%%%%%%%%%%%%%%%%%%%%%%%%%%%%%%%%%%%%%%%%%%%%%%%%%%%%%%%%%%%%%%%%%%%%%%%%%%%%%%%%%%%%%%%%%%%%%%%%%%%%%%%%%%%%%%%%%%%%%%%%%%%%%%%%%%%%%%%%%%%%%%%%%%%%%%%%%%%%%%%%%%%%%%%%%%%%%%%%%%%%%%%%%%%%%%%%%%%%%%%%%%%%%%%%%%%%%%%%%%%%%%%%%%%%%%%%%%%%%%

\section{\bf Proof of Theorem \ref{t91}}

  \noindent  \emph{Proof of $(i)$}. If $f^2 = \textrm{id}$, there is nothing to prove. 
 If $f^2 \neq id$, from Theorem \ref{t52}.(1), it suffices to find two involutions $\tau$ and $h$ in $\textrm{PL}^+(\mathbb{S}^1)$ 
such that Fix$(\tau f) \neq \emptyset$ and $\Delta_{\tau f} = -\Delta_{\tau f} \circ h$ since in that case, $f$ is a composition of three involutions of 
$\textrm{PL}^+(\mathbb{S}^1)$. There is a point $x$ in $\mathbb{S}^1$  
such that $x\neq f^2(x)$. We can assume that the points $x, \ f(x)$ and $f^2(x)$ occur in that order anticlockwise around 
$\mathbb{S}^1$. Choose a point $y$ in $(x, f(x))$ such that $f^{-1}(y)$ be in $(f^2(x), x)$.
 Let $u : [x, f(x)] \longrightarrow [f(x), x]$ be an orientation-preserving piecewise linear homeomorphism such that:\\
 $u(y) = f(y)$,
 $$\begin{array}{ll}
u(t) < \min(f(t), f^{-1}(t)), &  \textrm{for} \ t \in (x, y); \\
f(t)<u(t)<f^{-1}(t), & \textrm{for }
 \
 t\in (y, f(x)).
\end{array}$$
 \noindent Then, let $\tau$ be the involution in PL$^+(\mathbb{S}^1)$ defined by $$\tau(t) = \left\{
                                                                                 \begin{array}{ll}
                                                                                   u(t), &  \textrm{if } \ t \in [x, f(x)] \\
                                                                                    u^{-1}(t), & \textrm{if
                                                                                    }
                                                                                    \
                                                                                                              t \in [f(x), x]
                                                                                     \end{array}
                                                                                       \right.$$
   \noindent We have $\tau f(t) = t$ if and only if $t=x$ or $t=y$. So, $\textrm{Fix}(\tau f) = \{x,  y\}$. Moreover, we have:
   $$\begin{array}{ll}
\forall \ t \in (x,y), &  f(t) > u(t) \Longleftrightarrow \tau f(t) > t \\
 \forall \ t \in (y, f(x)), & f(t) < u(t) \Longleftrightarrow \tau f(t) < t \\
\forall \ t \in (f(x), x), & u\left(u^{-1}(t)\right) < f^{-1}\left(u^{-1}(t)\right) \Longleftrightarrow \tau f(t) < t.
                                                                                      \end{array}$$
  \noindent Therefore: $$\Delta_{\tau f}(t) = \left\{
\begin{array}{ll}
 0, &  \textrm{if } \ t = x, y \\
1, & \textrm{if }
\
t \in (x, y) \\
-1, & \textrm{if } \ t \in (y, x).                        
\end{array}
\right.$$
                                                                                       
  Now, let $v: [x,y] \longrightarrow [y,x]$ be any orientation-preserving piecewise linear homeomorphism, and let $h$ be 
the involution in $\textrm{PL}^+(\mathbb{S}^1)$ defined by  $$h(t) = \left\{
                                                                                 \begin{array}{ll}
                                                                                   v(t), &  \textrm{if } \ t \in [x, y] \\
                                                                                    v^{-1}(t), & \textrm{if
                                                                                    }
                                                                                    \
                                                                                                              t \in [y, x]
                                                                                     \end{array}
                                                                                       \right.$$
  
  \noindent It is easy to see that $h$ satisfies $\Delta_{\tau f} \circ h = - \Delta_{\tau f}$.
   We conclude that each member
  of $\textrm{PL}^+(\mathbb{S}^{1})$ can be expressed as a composite of three involutions of
  $\textrm{PL}^+(\mathbb{S}^{1})$.
  So, $\textrm{PL}^+(\mathbb{S}^{1}) =
  I_3(\textrm{PL}^+(\mathbb{S}^{1})) =
  R_2(\textrm{PL}^+(\mathbb{S}^{1}))$.
  There are elements in $\textrm{PL}^+(\mathbb{S}^{1})$ which are not strongly reversible in $\textrm{PL}^+(\mathbb{S}^{1})$; one can choose,
  for example, a homeomorphism $f \in \textrm{PL}^+(\mathbb{S}^{1})$ which is not an involution and
  with rotation number $\rho(f) = \frac{1}{2}$ (such map $f$ is not strongly reversible in
  $\textrm{PL}^+(\mathbb{S}^{1})$ by Lemma \ref{l53}).
  The fact that  $R_1(\textrm{PL}^+(\mathbb{S}^{1})) {\subset}_{\neq} I_2(\textrm{PL}(\mathbb{S}^{1}))$
  follows from Theorem \ref{t11}.
\medskip

 \textit{ Proof of $(ii)$}. If $f\in \textrm{PL}^+(\mathbb{S}^{1})$, then by (i), $f\in I_3(\textrm{PL}(\mathbb{S}^{1}))$.
 If $f\in \textrm{PL}^-(\mathbb{S}^{1})$, then
  $\textrm{Fix}(f)\neq \emptyset$. Let $a\in \textrm{Fix}(f)$ and let $I = \mathbb{S}^{1}\setminus \{a\}$. 
Then, there exists an open interval $\widehat{I}$ in $\mathbb{R}$ such that the map 
$\varphi: \widehat{I}\longrightarrow I$ given by $\varphi(t) = e^{2i\pi t}$, is a homeomorphism. Set $g = \varphi^{-1} f \varphi$. 
We have $g \in \textrm{PL}^-(\widehat{I})$. Choose an involution $\sigma \in \textrm{PL}^-( \widehat{I})$ such that, 
for each $x\in \widehat{I}$, $\sigma(x) > g(x)$. Then $g(\sigma(x)) < x$, for each $x\in \widehat{I}$. Therefore, $g\sigma$ is a 
fixed point free element in $\textrm{PL}^+( \widehat{I})$. By Lemma 4.2, (1), it is strongly reversible 
in $\textrm{PL}( \widehat{I})$ by an element of $\textrm{PL}^-( \widehat{I})$; which means that there exist three involutions $u$, $v$ in $\textrm{PL}^-(\widehat{I})$ such that 
$g\sigma = u v$. Thus $g = uv\sigma$. It follows that $f_{| I} = \varphi g \varphi^{-1} = 
 (\varphi u \varphi^{-1})(\varphi v \varphi^{-1})(\varphi \sigma \varphi^{-1})$. By extending $\varphi u \varphi^{-1}$,  
$\varphi v \varphi^{-1}$ and $\varphi \sigma \varphi^{-1}$ to $\mathbb{S}^{1}$ by fixing $a$, we get three involutions $\tau_1, \ \tau_2$ and $\tau_3$ 
in $\textrm{PL}^-(\mathbb{S}^{1})$ satisfying 
$f = \tau_1\tau_2\tau_3$. Hence $f\in I_3(\textrm{PL}(\mathbb{S}^{1}))$. We conclude that
$\textrm{PL}(\mathbb{S}^{1}) = I_3(\textrm{PL}(\mathbb{S}^{1}))=
R_2(\textrm{PL}(\mathbb{S}^{1}))$. Moreover,
$\textrm{PL}(\mathbb{S}^{1})\neq I_2(\textrm{PL}(\mathbb{S}^{1}))$
as in the proof of (i). From Theorem \ref{t2}, we have
$R_1(\textrm{PL}(\mathbb{S}^{1})) =
I_2(\textrm{PL}(\mathbb{S}^{1}))$. Now to show that
$I_2(\textrm{PL}(\mathbb{S}^{1})) \neq
I_1(\textrm{PL}(\mathbb{S}^{1}))$, it suffices to choose a
nontrivial reversible element $f$ in $\textrm{PL}^+(\mathbb{S}^{1})$
which is not the identity and with rotation number $\rho(f) = 0$.
\qed
\medskip

\begin{rem} \label{r4} Contrarily to $\mathrm{PL}^-(\mathbb{R)}$ (cf. Lemma \ref{l299}), there exists an element of $\mathrm{PL}^-(\mathbb{S}^{1})$ which is strongly
reversible in $\mathrm{PL}(\mathbb{S}^{1})$ but not an involution.
\end{rem}

\begin{proof} Indeed, Suppose that the remark is not true. We will prove in this case that any element 
$f\in \textrm{PL}^-(\mathbb{S}^{1})$ is an involution, this leeds to a contradiction since there are elements in
$\textrm{PL}^-(\mathbb{S}^{1})$ which are not involutions.
 Indeed, let $\sigma$ be any involution in $\textrm{PL}^-(\mathbb{S}^{1})$. Then $\sigma f\in \textrm{PL}^+(\mathbb{S}^{1})$ and from Theorem \ref{t91}, (i), there exist three involutions 
$\tau_1, \ \tau_2$ and $\tau_3$ in $\textrm{PL}^+(\mathbb{S}^{1})$ such that $\sigma f = \tau_1 \tau_2 \tau_3$. This implies that 
$f = \sigma \tau_1 \tau_2 \tau_3$. By assumption, $\sigma \tau_1$ is an involution in 
$\textrm{PL}^-(\mathbb{S}^{1})$ and then so is $(\sigma \tau_1) \tau_2$. We conclude that $f$ is an involution.

\end{proof}
%\textbf{Question}

%%% ----------------------------------------------------------------------

\end{document}